\documentclass[12pt,leqno]{article}
\usepackage{a4wide,epsf,latexsym}
\usepackage[centertags]{amsmath}
\usepackage{graphicx}
\numberwithin{equation}{section}

\newtheorem{remark}{Remark}

\newcommand{\ve}{\varepsilon}
\newcommand{\divergenz}{\ensuremath{\operatorname{div}} }
\newcommand{\diag}{\ensuremath{\operatorname{diag}} }

\begin{document}

\title {Some remarks on strongly coupled systems of convection-diffusion
equations in 2D}

\author{Hans-G. Roos, TU Dresden}

\date{13.2.2015}

\maketitle

\begin{abstract}
Almost nothing is known about the layer structure of solutions to strongly
coupled systems of convection-diffusion equations in two dimensions.
In some special cases we present first results.
\end{abstract}

{\it AMS subject classification}: 65 L10, 65 L12, 65 L50

\section{Introduction}
In the survey paper \cite{LS10} the authors present some results for
weakly coupled systems in one and two space dimensions but state
concerning strongly coupled convection-diffusion problems ,,we have only a limited
grasp of the situation''. We aim to provide an inside into the nature of
such problems at least in some special cases.

A practical example of strongly coupled systems of convection-diffusion
equations in 2D (related to magnetohydrodynamic duct flow)
 is numerically studied in \cite{HSY11}, namely
\begin{gather}\label{1.1a}
	\begin{aligned}
		-\varepsilon \Delta u+a\nabla b &=f_1,\\
		-\varepsilon \Delta b+a\nabla u &=f_2
	\end{aligned}
\end{gather}
with some boundary conditions.

Let us more generally consider the vector-valued function $u=(u_1,u_2)^T$ solving the system
\begin{subequations}\label{(1)}
  \begin{align}
-\ve \Delta u +A_1\frac{\partial u}{\partial x_1}+A_2\frac{\partial u}{\partial x_2}
   +\rho\, u&=f
\qquad \mbox{in}\,\, \Omega,\\
u&=0 \qquad \mbox{on}\,\, \Gamma= \partial \Omega,
   \end{align}
\end{subequations}
where $\ve$ is a small positive parameter. We assume the matrices
$A_1,A_2$ to be symmetric and $C^1$ and that the unit outer normal $\nu = (\nu_1,\nu_2)$ to $\Omega$ exists a.e.~on  $\partial\Omega$.

Because
\[
(\sum A_i\frac{\partial u}{\partial x_i},u)=
  \frac{1}{2}\int_{\Gamma}(\nu\cdot Au,u)d\Gamma-\frac{1}{2}\big((\divergenz A)u,u\big)
\]
with
\[
    \divergenz A=\sum \frac{\partial A_i}{\partial x_i},
\]
it is standard to assume
\begin{equation}
  \rho>  \frac{1}{2}\sup_{\Omega}\|\divergenz A\|_\infty.
\end{equation}
Then, problem \eqref{(1)} admits a unique weak solution.

To describe the reduced problem we introduce the matrix
\begin{equation}
  B:=\nu_1 A_1+\nu_2 A_2
\end{equation}
Suppose $B$ to be nonsingular,
i.e., $\partial\Omega$ to be noncharacteristic. Then $B$ allows the
decomposition
\begin{equation}\label{dec}
  B=B^{+} +B^{-} ,
\end{equation}
where $B^{+}$
 is positive semidefinite, $B^{-}$ negative semidefinite
 and the eigenvalues of $B^{+}$ are the positive eigenvalues of
 $B$ and $0$. The reduced
 problem to \eqref{(1)} is then given by
\begin{equation}
  A_1\frac{\partial u_0}{\partial x_1}+A_2\frac{\partial u_0}{\partial x_2}
   +\rho\, u_0=f
\qquad \mbox{with}\,\,
B^{-}u_0=0 \qquad \mbox{on}\,\, \Gamma.
\end{equation}
In \cite{Si70} it was proved that $u$ converges for $\varepsilon\to 0$ to
$u_0$. But concerning the convergence rate we only know the result
of  \cite{LS74} for a problem with different boundary conditions:
for $f\in H^1$ one has in the $L_2$ norm
\begin{equation}
  \|u-u_0\|_0\le C{\ve}^{1/2}\|f\|_1.
\end{equation}
In the literature not much is known about the structure of layers.
In the next Section we will discuss this question in, at least, some
special cases and consider problems with two small
parameters, as well.

\section{The reduced problem and location of layers}
Let us first revisit the situation in 1D studied in \cite{RR10}.
We take the example
\begin{equation}
-\ve u'' -	\begin{pmatrix}
			3 & 4\\
			4 & -3
		\end{pmatrix} u' =	\begin{pmatrix}
						1\\
						2
					\end{pmatrix}
\end{equation}
with $u(0)=u(1)=0$. In the one parameter case the eigenvalues and
eigenvectors of $A$ determine the asymptotic structure of the
solution. Here, the eigenvalues are 5 and -5 with the corresponding
eigenvectors $(2,1)^T$ and $(1,-2)^T$. With that knowledge and the
general solution of the reduced equation we made in  \cite{RR10}
the Ansatz
\begin{equation}\label{2.5}
u_{as}^0 = w_0 + w_h \begin{pmatrix} c_1\\c_2 \end{pmatrix} + d_1	 \begin{pmatrix}
													2\\
													1
												\end{pmatrix} \exp(-5 x/\ve)
										+ d_2	\begin{pmatrix}
													1\\
													-2
												\end{pmatrix} \exp\big(-5 (1-x)/\ve\big).
\end{equation}
The four constants $c_1,c_2,d_1,d_2$ are computed from the boundary
conditions. Eliminating all exponentially small terms one gets indirectly
the solution of the reduced problem, for our example
\[
  u_0=\begin{pmatrix}
	11/25\, x-8/25							\\
	-2/25\,x-4/25												
												\end{pmatrix}.
\]
Observe that $u_0$ does not satisfy either of the given boundary conditions.

There are two possibilities to derive the boundary conditions for
the reduced equation {\it directly}. First as discussed in the Introduction,
we can decompose the matrix $B$.

In our example we have
\[
(Av,v)=3v_1^2+8v_1v_2-3v_2^2=(2v_1+v_2)^2-(v_1-2v_2)^2.
\]
Therefore, the reduced problem is given by
\[ -	\begin{pmatrix}
				3 & 4\\
				4 & -3
			\end{pmatrix} u_0' =	\begin{pmatrix}
									1\\
									2
								\end{pmatrix}\]
with the boundary conditions
\begin{gather}\label{eq:boundaryconditions}
   (2u_0^1+u_0^2)(1)=0 \quad{\rm and}\quad
   (u_0^1-2u_0^2)(0)=0.
\end{gather}

An alternative approach consists in the diagonalization of the
matrix $A$. Set
\[
u=Tv \quad {\rm or}\quad v=T^{-1}u.
\]
After transformation $T^{-1}AT$ is a diagonal matrix, we obtain a decoupled
system in $v$. In our example,
\[
u=Tv\quad {\rm with}\quad T=\frac{1}{\sqrt 5}\begin{pmatrix}
				2 & 1\\
				1 & -2
			\end{pmatrix}
\]
generates
\begin{subequations}
  \begin{align}
-\ve v_1''-5v_1'
   &=\frac{4}{\sqrt 5},
\\
-\ve v_2''+5v_2'
   &=-\frac{3}{\sqrt 5}.
   \end{align}
\end{subequations}
>From the sign of the coefficients of the first order derivatives one
can conclude where the layers of $v$ are located. In our example
$v_1$ has a layer at $x=0$, $v_2$ at $x=1$. The back transformation
yields:
\begin{itemize}
\item the correct boundary conditions for the reduced problem are given by \eqref{eq:boundaryconditions}
\item both components $u_1$ and $u_2$ of the solution have layers at $x=0$ and $x=1$.
\end{itemize}

The transformation technique can also be applied for special systems
in two dimensions. Let us assume that the constant matrices
$A_1, A_2$ of the system \eqref{(1)} admit the representation
\begin{equation}
A_1=T\, \diag(\lambda_1^1,\lambda_1^2))\,T^{-1},\quad
A_2=T\, \diag(\lambda_2^1,\lambda_2^2))\,T^{-1}
\end{equation}
with the orthogonal matrix
\[
T=\begin{pmatrix}
	\sin \phi & \cos \phi						\\
	-\cos \phi & \sin \phi	\end{pmatrix}
\]
for some $\phi \in [0,2\pi]$. Remark that this is the case if and only if $A_1A_2=A_2A_1$.

Then, the transformation
\[
u=Tv \quad {\rm or}\quad v=T^{-1}u
\]
yields also a decoupled system, namely (for $\rho=0$, but this is not
essential)
\begin{subequations}
  \begin{align}
-\ve \Delta v_1 +\lambda_1^1\frac{\partial v_1}{\partial x_1}+
   \lambda_2^1\frac{\partial v_1}{\partial x_2}
   &=\tilde f_1\\
-\ve \Delta v_2 +\lambda_1^2\frac{\partial v_2}{\partial x_1}+
   \lambda_2^2\frac{\partial v_2}{\partial x_2}
   &=\tilde f_2.
   \end{align}
\end{subequations}
That means: the sign of $(\lambda_1^1,\lambda_2^1)^T\cdot \nu$ determines
the location of the layers of $v_1$, the sign of
$(\lambda_1^2,\lambda_2^2)^T\cdot \nu$ the location of the layers of
$v_2$. Let us introduce
\[
\Gamma_1^{+}=\{x\in\Gamma: \,\,(\lambda_1^1,\lambda_2^1)\cdot \nu>0
\}
\]
and analogously $\Gamma_1^{-},\Gamma_1^{0},\Gamma_2^{+},\Gamma_2^{-},
\Gamma_2^{0}$. Then, the conditions
\[
T^{-1}_{r1}u|_{\Gamma_1^{-}}=0\quad {\rm and}\quad
    T^{-1}_{r2}u|_{\Gamma_2^{-}}=0
\]
($T^{-1}_{r1}$ denotes the first row of $T^{-1}$)
are the correct boundary conditions
for the reduced equation.

\begin{figure}
	\centering%
	\includegraphics{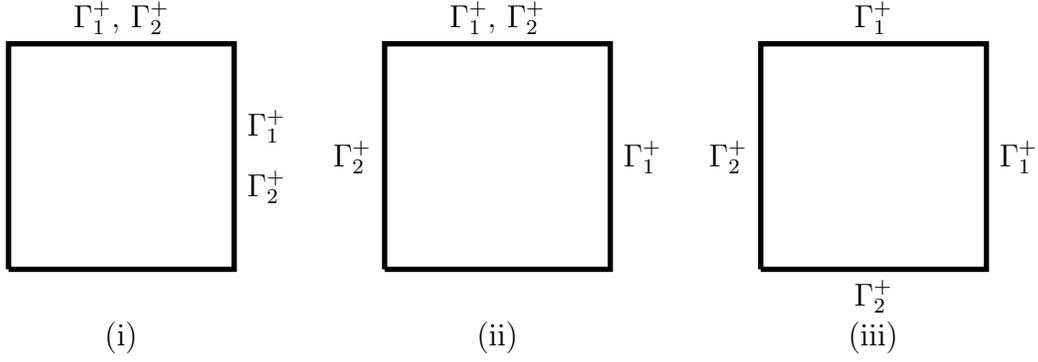}%
	\caption{Location of boundary layers}%
	\label{fig:layerlocation}
\end{figure}

Assume $\Omega=(0,1)^2$. Then we have three typical cases (see Figure \ref{fig:layerlocation}).
\begin{itemize}
\item (i) all $\lambda$ are positive\\
  Then $u_1$ and $u_2$ have overlapping layers at $x_1=1$ and at
  $x_2=1$.
\item (ii) only $\lambda_1^2$ is negative\\
  Then we have again overlapping layers at $x_2=1$, but ordinary layers
  of $u_1$ and $u_2$ at $x_1=0$ and at $x_1=1$.
\item (iii) both components of $\lambda_1$ and $\lambda_2$ have a
     different sign\\
     Then $u_1$ and $u_2$ have ordinary layers at every edge of the unit square.
\end{itemize}

Let us remark that the duct flow example mentioned in the introduction
under the assumption $a<0$ of \cite{HSY11} falls under case (iii) with
\[
T=\frac{1}{\sqrt 2}\begin{pmatrix}
	1\,\,\, 1						\\
	-1\,\,\,1											
												\end{pmatrix}.
\]

Let us now consider the two parameter case
\begin{subequations}\label{two}
  \begin{align}
-{\cal E} \Delta u +A_1\frac{\partial u}{\partial x_1}+A_2\frac{\partial u}{\partial x_2}
   +\rho\, u&=f
\qquad \mbox{in}\,\, \Omega\\
u&=0 \qquad \mbox{on}\,\, \Gamma= \partial \Omega
   \end{align}
\end{subequations}
with ${\cal E}=\diag(\ve_1,\ve_2)$ and two small positive parameters. We are
mostly interested in the case $\ve_1/\ve_2 \ll 1$.

Let us first collect some facts from \cite{Ro12} for the one-dimensional
system
\begin{equation}\label{two1D}
 - \diag(\varepsilon_1,\varepsilon_2)u''-Bu'+Au=f,\quad u(0)=u(1)=0.
\end{equation}
Assume
\[
B=\begin{pmatrix}  -b_{11} & b_{12} \\ b_{12} & b_{22}  \end{pmatrix}\quad\mbox{with}\quad  b_{11}>0 \quad {\rm and}\,\,
  \hat D=b_{11}b_{22}+b_{12}^2>0.
\]
Surprisingly, not the eigenvalues of $B$ generate the structure of
the boundary layers. The layer at $x=0$ is characterized by the
exponential $\exp(-b_{11}/\ve_1)$, but the layer at $x=1$ by the
exponential $\exp(\frac{\hat D}{b_{11}}\frac{1}{\ve_2})$.

The quantities $-b_{11}$ and $\frac{-\hat D}{-b_{11}}$ are quotients
of the leading principal minors of the matrix $B$. These values arise as well
in the $LDL^T$ decomposition of the matrix $B$. In accordance with
that observation the authors of \cite{MOPX11} proposed to introduce
new variables $v$ by $v=L^Tu$ resulting in
\begin{equation}
 - {\cal E^*}v''-\diag(d_1,d_2)v'+A^*v=f^*,\quad v(0)=v(1)=0
\end{equation}
with
\[
  {\cal E^*}=-L^{-1}{\cal E}(L^{-1})^T.
\]
Unfortunately, the matrix ${\cal E^*}$ is not diagonal. But its structure
\[
{\cal E^*}=\begin{pmatrix}  \ve_1 & -\ve_1\,l \\
 -\ve_1\,l & \ve_1\,l^2+\ve_2
 \end{pmatrix}
\]
tells us (remember  $\ve_1/\ve_2\ll1$) the following:
Assuming, for instance, $d_1<0$ and $d_2>0$, both solution components
have strong layers at $x=0$ of the structure $\exp\big((d_1/\ve_1)x\big)$. However,
at $x=1$ a {\it weak} layer of $v_2=u_2$ of the structure
 $O(\ve_1/\ve_2)\exp\big(-(d_2/\ve_2)(1-x)\big)$ generates a strong layer of $v_1$
and, consequently, $u_1$. In \cite{Ro12} this solution behavior was derived
from some asymptotic approximation.

The signs of $d_1,d_2$ also define the boundary conditions of the
reduced problem, under our assumptions
\begin{equation}
u_0^2(1)=0,\quad {\rm and}\quad u_0^1(0)+lu_0^2(0)=0.
\end{equation}
Consider the example
\[
 - \diag(\varepsilon_1,\varepsilon_2)u''-\begin{pmatrix}
				3 & 4\\
				4 & -3
			\end{pmatrix} u' =	\begin{pmatrix}
									1\\
									2
								\end{pmatrix}
 \]
with $u(0)=u(1)=0$.
Here
\[
L=\begin{pmatrix}
				1 & 0\\
				4/3 & 1
			\end{pmatrix}
\]
and $D=\diag(3,-25/3)$. Consequently, in the case $\ve_1\ll\ve_2$
we conclude
\begin{itemize}
\item the boundary conditions for the reduced problem are
  $$u_0^2(1)=0,\quad {\rm and}\quad u_0^1(0)+4/3u_0^2(0)=0
  $$
\item the solution u is characterized by the layer terms
$$\begin{pmatrix}
									4\\
									-3
								\end{pmatrix}\exp(-\frac{3}{\ve_1}x)
\quad {\rm and}\quad
\begin{pmatrix}
									3\\
									O(\ve_1/\ve_2)
								 \end{pmatrix}\exp\Big(-\frac{25}{3\ve_2}(1-x)\Big).$$
\end{itemize}

\begin{remark} In the one parameter case the decomposition \eqref{dec}
allows to define the boundary conditions for the reduced problem. But
with two parameters, for instance in 1D, the matrix ${\cal E}^{-1}B$
is not symmetric. Symmetrization of that matrix with a transformation
of the type
\[
T^{-1}({\cal E}^{-1}B)T \quad \text{with}\quad T = \begin{pmatrix}
				1 & \mu\\
			(\frac{\ve_1}{\ve_2})^{1/2}\mu & 1
			\end{pmatrix}
\]
and adequately chosen $\mu$ is possible. But it is simpler to
apply the $LDL^T$ decomposition described above.
$\Box$
\end{remark}

The results for the two parameter problem \eqref{two1D} can now be
generalized to \eqref{two} if we assume
\begin{equation}
A_1=L\, \diag(d_1^1,d_1^2)\,L^T,\quad
A_2=L\, \diag(d_2^1,d_2^2)\,L^T.
\end{equation}
The transformation
\[
u=(L^T)^{-1}v \quad {\rm or}\quad v=L^Tu
\]
yields the system
\begin{gather}
	\begin{aligned}
-\ve_1 \Delta v_1-\ve_1 q\Delta v_2+d^1\cdot\nabla v_1+\tilde \rho(v_1+\alpha v_2)&=\tilde f_1,\\
-\ve_1 q\Delta v_1-(\ve_1q^2+\ve_2)\Delta v_2
+d^2\cdot\nabla v_2+ \rho  v_2&=\tilde f_2.
	\end{aligned}
\end{gather}
Next we introduce as above
$\Gamma_1^{-},\Gamma_1^{+},\Gamma_2^{-},\Gamma_2^{+}
$ with respect to the vectors $d^1,d^2$. Both components $v_1$ and $v_2$ as well as
$u_1$ and $u_2$ have strong exponential layers related to $\ve_1$ on $\Gamma_1^{+}$.
But on $\Gamma_2^{+}$ the component $v_2=u_2$ has only a weak layer.
That weak layer generates a strong layer of $u_1$ related to the
larger parameter $\ve_2$.

The correct boundary conditions for the reduced equation are
\[
L^T_{r1}u|_{\Gamma_1^{-}}=0\quad {\rm and}\quad
    u_2|_{\Gamma_2^{-}}=0.
\]

Summarizing: Under some restrictive conditions it is possible to characterize the
possible layers of strongly coupled convection-diffusion systems in 2D. Moreover,
the adequate boundary conditions for the reduced problem can be derived (important
for numerical methods which first solve the simpler reduced problem).

\end{document}